\documentclass[conference]{IEEEtran}
\ifCLASSINFOpdf
\else
\fi
\usepackage{amsmath,amssymb}
\newtheorem{definition}{Definition}
\newtheorem{theorem}{Theorem}

\begin{document}
\title{Numerical Solution of Fractional Control Problems via Fractional Differential Transformation}

\author{\IEEEauthorblockN{Josef Rebenda}
\IEEEauthorblockA{CEITEC BUT\\
Brno University of Technology\\
Purky\v{n}ova 123\\
612 00 Brno, Czech Republic\\
Email: josef.rebenda@ceitec.vutbr.cz}
\and
\IEEEauthorblockN{Zden\v{e}k \v{S}marda}
\IEEEauthorblockA{CEITEC BUT\\
Brno University of Technology\\
Purky\v{n}ova 123\\
612 00 Brno, Czech Republic\\
Email: smarda@feec.vutbr.cz}
}


%Telephone: (800) 555--1212, Fax: (888) 555--1212}

\IEEEspecialpapernotice{978-0-7695-6213-1/17 \copyright\ 2017 IEEE. DOI 10.1109/EECS.2017.29 Personal use of this material is permitted. Permission from IEEE must be obtained for all other uses, in any current or future media, including reprinting/republishing this material for advertising or promotional purposes, creating new collective works, for resale or redistribution to servers or lists, or reuse of any copyrighted component of this work in other works.}

\maketitle

\begin{abstract}
In the paper we deal with linear fractional control problems with constant delays in the state. Single-order systems with fractional derivative in Caputo sense of orders between 0 and 1 are considered. The aim is to introduce a new algorithm convenient for numerical approximation of a solution of the studied problem. The method consists of the fractional differential transformation in combination with general methods of steps. The original system is transformed to a system of recurrence relations. Approximation of the solution is given in the form of truncated fractional power series. The choice of order of the fractional power series is discussed and the order is determined in relation to the order of the system. An application on a two-dimensional fractional system is shown. Exact solution is found for the first two intervals of the method of steps. The result for Caputo derivative of order 1 coincides with the solution of first-order system with classical derivative. We conclude that the algorithm is applicable, efficient and gives reliable results.
\end{abstract}


% \fi
\IEEEpeerreviewmaketitle

\section{Introduction}
Fractional-order derivatives are a generalization of integer-order derivative. Different fractional derivatives have been
defined in fractional calculus, which is studied in detail in monographs of
Oldham and Spanier \cite{oldham}, Miller and Ross \cite{miller},
Samko et al. \cite{samko}  or Das \cite{das}.

Mathematical modeling of
systems and processes with the use of fractional
derivatives leads to fractional differential equations. Theory and applications of fractional differential equations are covered by monographs of Podlubny \cite{podlubny}, Kilbas et al.
\cite{kilbas} and Diethelm \cite{diethelm}.
Fractional differential equations and systems occur in mathematical
models of mechanical, biological,
chemical, physical and medical phenomena as well as in other areas of real life. It has become apparent
that fractional-order models reflect the behavior of many
real-life processes more accurately than integer-order
ones. For more details concerning fractional calculus
and its practical applications we refer to the monographs mentioned above.

Fractional systems with delays in the
state were discussed by Sikora \cite{sikora2003,sikora2016} and Buslowicz
\cite{buslowicz}. Fractional systems with delays
in control were analyzed by Sikora \cite{sikora2005}, Trzasko \cite{trzasko}, Kaczorek
\cite{kaczorek2011}, Balachandran et al. \cite{balachandran2012a,balachandran2012c} as well as Kaczorek and Rogowski \cite{kaczorek2015}.

In some cases, often in applications, it is not possible to find a solution of dynamical systems with delays in state or control in analytical form. Therefore, it is important to find a way how to approximate solutions of such systems numerically. Many convenient methods can be found in literature, e.g. in the monographs of Bellen and Zennaro \cite{bellen}, Sun and Ding \cite{sun}, or recent papers by Guglielmi and Hairer \cite{guglielmi} and Rebenda and \v{S}marda \cite{rebenda2017a}.

In the last years, numerous research papers are dedicated to study of numerical methods of solving fractional differential equations and related problems. Indicatively we mention papers by Jannelli et al. \cite{jannelli}, Yang et al. \cite{yang} and Odibat et al. \cite{odibat}.

Specifically, publications about the differential transformation are growing in both number and quality during the last years. To indicate the progress in the field we mention recent papers by \v{S}amajov\'{a} and Li \cite{samajova}, \v{S}marda and Khan \cite{smarda}, Yu et al. \cite{yu}, Rebenda and \v{S}marda \cite{rebenda2017b} and Rebenda et al. \cite{rebenda2017c}.

Theory on the method of steps can be found for instance in monographs of Hale and Verduyn Lunel \cite{hale} or Kolmanovskii and Myshkis \cite{kolmanovski})

The purpose of this paper is to introduce a new algorithm for numerical solution of linear fractional optimal control problems with multiple time-delays in the state functions. The algorithm is a combination of the method of steps and the fractional differential transformation (FDT). The goal is to derive results that are applicable to concrete problems in real world.

The paper is organized as follows. First we define Caputo fractional derivative and introduce the linear fractional control problem in Section \ref{problem}. We continue with definition and properties of the fractional differential transformation in Section \ref{fdt}. Section \ref{mos} contains a method how to eliminate delays from the state provided the delays are commensurate. Then in Section \ref{order} we discuss the choice of the order of the truncated fractional power series used to approximate the solution. An illustratory example is presented in Section \ref{applications}. Finally, conclusions are made in Section \ref{conclusion}.

%%%%%%%%%%%%%%%%%%%%%%%%%%%%%%%%%%%%%%%%%%%%%%%%%%%%%%%

\section{Problem Statement}\label{problem}
In this paper we work with the
Caputo fractional derivative. The purpose is to avoid fractional initial conditions and to use integer-order initial conditions which have clear practical meaning.
\begin{definition}
The fractional derivative in Caputo sense is defined by
\begin{equation}\label{7}
_{t_0}^{C} \! D_t^{\alpha} f(t) = \frac{1}{\Gamma(n-\alpha)} \int_{t_0}^t \frac{f^{(n)}(s)}{(t-s)^{1+\alpha-n}} ds.
\end{equation}
where $n-1 \leq \alpha <n$, $n \in \mathbb{N}$, $t >t_0$.
\end{definition}

\subsection{Considered Problem}
We consider linear fractional control problem with constant delays in the state
\begin{equation}\label{1}
_0^{C} \! \mathbf{D}_t^{\nu} \mathbf{x} (t) = \mathbf{A}_0 (t) \mathbf{x}(t) + \sum\limits_{i=1}^{r} \mathbf{A}_i (t) \mathbf{x} (t - \tau_i) + \mathbf{B}(t) \mathbf{u}(t),
\end{equation}
for $t \geq 0$, where
\begin{itemize}
\item
$_0^{C} \! \mathbf{D}_t^{\nu} \mathbf{x} (t)$ denotes the $n$-dimensional vector of fractional derivatives of order $\nu$ in the Caputo sense, i.e. $_0^{C} \! \mathbf{D}_t^{\nu} \mathbf{x} (t) = \Bigl( {_0^{C} \! D}_t^{\nu} x_1 (t), \ldots, {_0^{C} \! D}_t^{\nu} x_n (t) \Bigr)^T$,
\item
$0 < \nu \leq 1$,
\item
$\mathbf{x}(t) \in \mathbb{R}^n$ is the $n$-dimensional state vector,
\item
$\mathbf{A}_i (t)$ are $n \times n$ matrices of real functions, $i=0, 1, \dots, r$,
\item
$\tau_i$ are constant delays in the state, $i=0, 1, \dots, r$,
\item
$\mathbf{B}(t)$ is an $n \times m$ matrix of real functions,
\item
$\mathbf{u}(t)$ is the $m$-dimensional vector of control functions.
\end{itemize}
Let $\tau = \max \{\tau_1,\tau_2,\dots,\tau_r\}$. A vector function $\mathbf{\Phi} (t) = (\phi_1 (t), \ldots, \phi_n (t) )^T$ needs to be assigned to system \eqref{1} on the interval $[-\tau, 0]$. This function is called initial complete state of the fractional differential system \eqref{1}. Furthermore, for the sake of simplicity, we assume that $\phi_j (t) \in C^1 ([-\tau,0])$ for $j=1, \dots, n$.

%%%%%%%%%%%%%%%%%%%%%%%%%%%%%%%%%%%%%%%%%%%%%%%%%%%%%

\section{Fractional Differential Transformation}\label{fdt}

We introduce the fractional differential transformation (FDT) in this section.
\begin{definition}\label{def1}
Fractional differential transformation of order $\alpha$ of a real function $u(t)$ at a point $t_0 \in \mathbb R$ in Caputo sense is $^{C} \! \mathcal{D}_{\alpha} \{ u(t) \} [t_0] = \{ U_{\alpha} (k) \}_{k=0}^{\infty}$,
where $k \in \mathbb{N}_0$ and $U_{\alpha} (k)$, the fractional differential transformation of order $\alpha$ of the $(\alpha k)$th derivative of function $u(t)$ at $t_0$, is defined as
\begin{equation}\label{c2}
U_{\alpha} (k) = \frac{1}{\Gamma (\alpha k + 1)} \left[ _{t_0}^{C} \! D_t^{\alpha k} u(t) \right]_{t=t_0},
\end{equation}
provided that the original function $u(t)$ is analytical in some right neighborhood of $t_0$.
\end{definition}
\begin{definition}\label{def2}
Inverse fractional differential transformation of $\{ U_{\alpha} (k) \}_{k=0}^{\infty}$ is defined using a fractional power series as follows:
\begin{equation}\label{c3}
u(t) = {^{C}} \! \mathcal{D}_{\alpha}^{-1} \Bigl\{ \{ U_{\alpha} (k) \}_{k=0}^{\infty} \Bigr\} [t_0]= \sum_{k=0}^{\infty} U_{\alpha} (k) (t-t_0)^{\alpha k}.
\end{equation}
\end{definition}

Convergence of the fractional power series \eqref{c3} in the definition of the inverse FDT is studied in \cite{odibat}. In applications, we will use some basic FDT formulas also listed in \cite{odibat}:
\begin{theorem} \label{t1}
Assume that $\{ F_{\alpha} (k) \}_{k=0}^{\infty}$, $\{ G_{\alpha} (k) \}_{k=0}^{\infty}$ and $\{ H_{\alpha} (k) \}_{k=0}^{\infty}$ are differential transformations of order $\alpha$ at $t_0$ of functions $f(t)$, $g(t)$ and
$h(t)$, respectively, and $r, \beta >0$.
\begin{align*}
\text{ If } f(t) &=  (t-t_0)^{r}, \text{ then } F_{\alpha} (k) =\delta \left( k - \frac{r}{\alpha} \right), \text{ where } \delta  
\text{ is}\\
\ & \qquad \qquad \qquad \text{ the Kronecker delta}. \\
\text{ If } f(t) &= g(t)h(t), \text{ then } F_{\alpha} (k) = \sum_{l=0}^k G_{\alpha} (l) H_{\alpha} (k-l).\\
\text{ If } f(t) &= \frac{g(t)}{(t-t_0)^r}, \text{ then } F_{\alpha} (k) = G_{\alpha} \left( k+ \frac{r}{\alpha} \right), \text{provided }\\
\ & \qquad \qquad \qquad G_{\alpha} (l) = 0 \text{ for } l < \frac{r}{\alpha},\\
\text{ If } f(t) &= {\displaystyle _{t_0}^{C} \! D_t^{\beta} g(t)}, \text{ then }\\
\ & \qquad \qquad \quad \ \ F_{\alpha} (k) = {\displaystyle \frac{\Gamma(\alpha k + \beta +1)}{\Gamma (\alpha k +1)}} G_{\alpha} \left( k+ \frac{\beta}{\alpha} \right). 
\end{align*}
\end{theorem}

%%%%%%%%%%%%%%%%%%%%%%%%%%%%%%%%%%%%%%%%%%%%%%%%%%%%%

\section{Elimination Of The Delays From The State}\label{mos}

In this section, we apply the method of steps to single-order fractional differential system with constant delays in the state \eqref{1}. For the rest of the paper, we assume that the constant delays are commensurate, which means that $\displaystyle{\frac{\tau_i}{\tau_j}}$ is a rational number for all pairs $i,j=1, \ldots, r$. We can assume that $0<\tau_1 < \ldots < \tau_r$. Now we define $k_i = \displaystyle{\frac{\tau_i}{\tau_1}}$ for $i = 1, \ldots, r$ and let $k_*$ be the least common multiple of denominators of the numbers $k_1, \ldots, k_r$. Finally, we set $\tau_* = \displaystyle{\frac{\tau_1}{k_*}}$. This $\tau_*$ is represents the length of the intervals on which we will look for approximate solutions in sequential steps.

To find the unique solution in the first interval, i.e. interval $(0,\tau_*]$, we substitute the initial complete state $\mathbf{\Phi}(t)$ in all places where state vector with delays appears. Then system \eqref{1} changes to a system of fractional differential equations without delays in the state
\begin{equation}\label{2}
_0^{C} \! \mathbf{D}_t^{\nu} \mathbf{x}_1 (t) = \mathbf{A}_0 (t) \mathbf{x}_1 (t) + \sum\limits_{i=1}^{r} \mathbf{A}_i (t) \mathbf{\Phi} (t - \tau_i) + \mathbf{B}(t) \mathbf{u}(t),
\end{equation}
subject to state vector at time $t=0$
\begin{equation}\label{3}
\mathbf{x}_1 (0) = \mathbf{\Phi} (0) = (\phi_1 (0), \ldots, \phi_n (0) )^T.
\end{equation}
In the second step, we find the solution $\mathbf{x}_2$ on interval $(\tau_*, 2 \tau_*]$. To the places in system \eqref{1} where state vector with delays occurs, we substitute either initial complete state $\mathbf{\Phi}(t)$ (if $t-\tau_i \in [-\tau,0]$ for $t \in (\tau_*, 2 \tau_*]$) or the state vector $\mathbf{x}_1 (t)$ computed in the first step (if $t-\tau_i \in (0,\tau_*]$ for $t \in (\tau_*, 2 \tau_*]$). System \eqref{1} can then either be the same as in \eqref{2} with the unknown state $\mathbf{x}_2$ or it can have the form
\begin{align}\label{4}
_0^{C} \! \mathbf{D}_t^{\nu} \mathbf{x}_2 (t) &= \mathbf{A}_0 (t) \mathbf{x}_2 (t) + \mathbf{A}_1 (t) \mathbf{x}_1 (t-\tau_1) \notag \\
 &+ \sum\limits_{i=2}^{r} \mathbf{A}_i (t) \mathbf{\Phi} (t - \tau_i) + \mathbf{B}(t) \mathbf{u}(t)
\end{align}
in interval $(\tau_*, 2 \tau_*]$ with state vector at time $t=\tau_*$
\begin{equation}\label{5}
\mathbf{x}_2 (\tau_*) = \mathbf{x}_1 (\tau_*) = (x_{11} (\tau_*), \ldots, x_{1n} (\tau_*) )^T.
\end{equation}
We continue with calculating solutions on further intervals in the same way. In the $j$th step, to the places with state vector with delay we substitute either the initial complete state $\mathbf{\Phi}(t)$ (if $t-\tau_i \in [-\tau,0]$ for $t \in ((j-1) \tau_*, j \tau_*]$) or the state vector $\mathbf{x}_l (t)$ computed in one of the previous steps (if $t-\tau_i \in ((l-1) \tau_*,l \tau_*]$ for $t \in ((j-1) \tau_*, j \tau_*]$).

%%%%%%%%%%%%%%%%%%%%%%%%%%%%%%%%%%%%%%%%%%%%%%%%%%%%%

\section{Determining The Order Of The Fractional Power Series}\label{order}

Without loss of generality, we demonstrate finding of optimal order of the FDT on single-order fractional differential system without delays in the state given by \eqref{2} subject to state vector at time $t=0$ given by \eqref{3}. Applying the FDT, in particular the formulas of Theorem \ref{t1}, to the system \eqref{2}, we obtain the following system of recurrence relations
\begin{align}\label{6}
{\displaystyle \frac{\Gamma(\alpha k + \nu +1)}{\Gamma (\alpha k +1)}} \mathbf{X}_{1 \alpha} \left( k+ \frac{\nu}{\alpha} \right) = &\sum\limits_{l=0}^k \mathbf{\mathcal{A}}_{0 \alpha} (l) \mathbf{X}_{1 \alpha} (k-l) \notag \\
 + \sum\limits_{i=1}^{r} \sum\limits_{l=0}^k \mathbf{\mathcal{A}}_{i \alpha} (l) \mathbf{F}_{i \alpha} (k-l) + &\sum\limits_{l=0}^k \mathbf{\mathcal{B}_{\alpha}}(l) \mathbf{U}_{\alpha}(k-l),
\end{align}
where $\mathbf{X}_{1 \alpha}$, $\mathbf{F}_{i \alpha}$ and $\mathbf{U}_{\alpha}$, respective $\mathbf{\mathcal{A}}_{0 \alpha}$, $\mathbf{\mathcal{A}}_{i \alpha}$ and $\mathcal{B}_{\alpha}$ are vectors, respective matrices of fractional differential transformations of order $\alpha$ at $0$.

Before we proceed with transformation of the initial conditions given by \eqref{3}, we need to determine the order of the fractional power series $\alpha$. For this purpose, we recall that our single-order system contains only fractional derivatives of order $\nu$.

From now on, we assume that $\nu \in \mathbb{Q}_+$. We can express $\nu$ as a fraction $\displaystyle{\frac{p}{q}}$ for some $p,q \in \mathbb{N}$. We look for $\alpha$ which satisfies the following conditions:
\begin{enumerate}
 \item
  $0 < \alpha \leq 1$.
 \item
  There is $k_{\nu} \in \mathbb{N}$ such that $\alpha \cdot k_{\nu} = \nu$.
\item
  There is $k_1 \in \mathbb{N}$ such that $\alpha \cdot k_1 = 1$.
\end{enumerate}
The last condition allows us to use polynomials as control functions.

There are infinitely many possibilities for the choice of $\alpha$. However, we propose that $\alpha$ should be chosen as largest as possible, which, in our case, is $\displaystyle{\frac{1}{q}}$ (reciprocal of the denominator of $\nu$).

The fractional differential transformation of the state vector at time $t=0$ given by \eqref{3} is then defined as
\begin{equation}\label{3'}
\mathbf{X}_{1 \alpha} (k)= \left\{ \begin{array}{ll}
\frac{1}{\Gamma (\alpha k +1)}  \left[ \frac{d^{\alpha k} \mathbf{x}_1 (t)}{dt^{\alpha k}} \right]_{t=0}, &  {\rm  if} \ \alpha k \in \mathbb{N}, \\
0, & {\rm  if} \ \alpha k \not \in \mathbb{N},
\end{array} \right.
\end{equation}
where $k= 0,1,2,\dots, (\frac{\lambda}{\alpha} -1)$ and $\lambda$ is the highest order of the considered fractional differential system, in our case $\lambda = \nu$. In particular, initial conditions \eqref{3} give us $\mathbf{X}_{1 \alpha} (0)= \mathbf{\Phi} (0)$.

The same procedure is applied to find the FDT of the system and the state vector in further intervals.

%%%%%%%%%%%%%%%%%%%%%%%%%%%%%%%%%%%%%%%%%%%%%%%%%%

\section{Applications}\label{applications}

We apply the algorithm to two-dimensional system studied in the paper by Rahimkhani et al. \cite{rahimkhani}
\begin{align}\label{ex1}
\begin{pmatrix}
_0^{C} \! {D}_t^{\nu} x_{1} (t)\\
_0^{C} \! {D}_t^{\nu} x_{2} (t)
\end{pmatrix}
&=
\begin{pmatrix}
t & 1\\
t & 2t
\end{pmatrix}
\begin{pmatrix}
x_1 \left( t-\frac{1}{3} \right)\\
x_2 \left( t-\frac{1}{3} \right)
\end{pmatrix} \notag \\
&+
\begin{pmatrix}
2 & t\\
t^2 & 0
\end{pmatrix}
\begin{pmatrix}
x_1 \left( t-\frac{2}{3} \right)\\
x_2 \left( t-\frac{2}{3} \right)
\end{pmatrix}
+
\begin{pmatrix}
0\\
1
\end{pmatrix}
u(t)
\end{align}
with initial complete states
\begin{equation}
x_1 (t) = x_2 (t) = 0 \quad \text{ for } t \in \left[ -\frac{2}{3}, 0 \right]
\end{equation}
and polynomial control function
\begin{equation}
u(t) = 2t +1.
\end{equation}
The fractional derivative $\nu$ can be an arbitrary rational number $\displaystyle{0 < \frac{p}{q} \leq 1}$.

First we rewrite system \eqref{ex1} as two equations which is more convenient for application of the procedure.
\begin{align}
_0^{C} \! {D}_t^{\nu} x_{1} (t) = t x_1 \left( t-\frac{1}{3} \right) &+ x_2 \left( t-\frac{1}{3} \right) \notag \\
+ 2 x_1 &\left( t-\frac{2}{3} \right) + t x_2 \left( t-\frac{2}{3} \right), \label{ex3}\\
_0^{C} \! {D}_t^{\nu} x_{2} (t) = t x_1 \left( t-\frac{1}{3} \right) &+ 2t x_2 \left( t-\frac{1}{3} \right) \notag \\
&+ t^2 x_1 \left( t-\frac{2}{3} \right) + 2t +1. \label{ex4}
\end{align}
Now we eliminate the delays from the state on the first interval $\displaystyle{\left[ 0, \frac{1}{3} \right]}$. On this interval, the system is very simple:
\begin{align*}
_0^{C} \! {D}_t^{\nu} x_{1} (t) = 0, \\
_0^{C} \! {D}_t^{\nu} x_{2} (t) = 2t +1.
\end{align*}
After performing FDT of order $\displaystyle{\alpha = \frac{1}{q}}$ at $t=0$ we get the system of recurrence relations
\begin{align}
X_{1, \alpha} (k+p) &= 0, \label{ex5} \\
X_{2, \alpha} (k+p) &= \frac{\Gamma \left( \frac{k}{q} +1 \right)}{\Gamma \left( \frac{k+p}{q} +1 \right)} \left( 2 \delta (k-q) + \delta (k) \right). \label{ex2}
\end{align}
The state vector at $t=0$ is identically zero, hence $X_{1, \alpha} (0) = 0 = X_{2, \alpha} (0)$. As we can see in \eqref{ex5} and \eqref{ex2}, the only nonzero coefficients of the fractional power series approximating the solution on $\displaystyle{\left[ 0, \frac{1}{3} \right]}$ are
\begin{align*}
X_{2, \alpha} (p) &= \frac{1}{\Gamma \left( \frac{q+p}{q} \right)} \quad (k=0), \\
X_{2, \alpha} (p+q) &= \frac{2}{\Gamma \left( \frac{2q+p}{q} \right)} \quad (k=q).
\end{align*}
In this case, the approximate solution coincides with the exact solution
\begin{align*}
x_{1, 1} (t) &= 0, \\
x_{2, 1} (t) &= \frac{1}{\Gamma \left( \frac{q+p}{q} \right)} t^{\nu} + \frac{2}{\Gamma \left( \frac{2q+p}{q} \right)} t^{\nu+1}.
\end{align*}
Following the steps of the procedure, we eliminate the delays from the state on the second interval $\displaystyle{\left[ \frac{1}{3}, \frac{2}{3} \right]}$. The system \eqref{ex3}, \eqref{ex4} changes to
\begin{align*}
_0^{C} \! {D}_t^{\nu} x_{1} (t) &= \frac{1}{\Gamma \left( \frac{q+p}{q} \right)} \left( t-\frac{1}{3} \right)^{\nu} + \frac{2}{\Gamma \left( \frac{2q+p}{q} \right)} \left( t-\frac{1}{3} \right)^{\nu+1}, \\
_0^{C} \! {D}_t^{\nu} x_{2} (t) &=  2t \left( \frac{1}{\Gamma \left( \frac{q+p}{q} \right)} \left( t-\frac{1}{3} \right)^{\nu} \right. \\
&\left. + \frac{2}{\Gamma \left( \frac{2q+p}{q} \right)} \left( t-\frac{1}{3} \right)^{\nu+1} \right) + 2t +1.
\end{align*}
FDT of order $\alpha$ but this time at $\displaystyle{t=\frac{1}{3}}$ leads to
\begin{align*}
X_{1, \alpha} (k&+p) = \frac{\Gamma \left( \frac{k}{q} +1 \right)}{\Gamma \left( \frac{k+p}{q} +1 \right)} \left( \frac{1}{\Gamma \left( \frac{q+p}{q} \right)} \delta (k-p) \right. \\
&+ \left. \frac{2}{\Gamma \left( \frac{2q+p}{q} \right)} \delta (k- (q+p)) \right), \\
X_{2, \alpha} (k&+p) = \frac{\Gamma \left( \frac{k}{q} +1 \right)}{\Gamma \left( \frac{k+p}{q} +1 \right)} \left( \frac{2}{3} \frac{1}{\Gamma \left( \frac{q+p}{q} \right)} \delta (k-p) \right. \\
+ &\left[ \frac{2}{\Gamma \left( \frac{q+p}{q} \right)} + \frac{2}{3} \frac{2}{\Gamma \left( \frac{2q+p}{q} \right)} \right] \delta \bigl( k- (q+p) \bigr) \\
+ & \left. \frac{4}{\Gamma \left( \frac{2q+p}{q} \right)} \delta \bigl( k- (2q+p) \bigr) + 2 \delta (k-q) + \frac{5}{3} \delta (k) \right).
\end{align*}
The state vector at $\displaystyle{t=\frac{1}{3}}$ is $(0, \frac{4}{9})^T$. The nonzero coefficients of the power series then are:
\begin{align*}
A &= X_{1, \alpha} (2p) = \frac{1}{\Gamma \left( \frac{q+2p}{q} \right)}, \\
B &= X_{1, \alpha} (2p+q) = \frac{2}{\Gamma \left( \frac{2q+2p}{q} \right)}, \\
C &= X_{2, \alpha} (0) = \frac{4}{9}, \\
D &= X_{2, \alpha} (p) = \frac{5}{3}\frac{1}{\Gamma \left( \frac{q+p}{q} \right)}, \\
E &= X_{2, \alpha} (p+q) = \frac{2}{\Gamma \left( \frac{2q+p}{q} \right)}, \\
F &= X_{2, \alpha} (2p) = \frac{2}{3}\frac{1}{\Gamma \left( \frac{q+2p}{q} \right)},
\end{align*}
\begin{align*}
G &= \! X_{2, \alpha} (2p+q) \! = \frac{{\Gamma \left( \frac{2q+p}{q} \right)}}{\Gamma \left( \frac{2q+2p}{q} \right)} \! \left[ \frac{2}{\Gamma \left( \frac{q+p}{q} \right)} + \frac{2}{3} \frac{2}{\Gamma \left( \frac{2q+p}{q} \right)} \right] \! , \\
H &= X_{2, \alpha} (2p+2q) = \frac{{\Gamma \left( \frac{3q+p}{q} \right)}}{\Gamma \left( \frac{3q+2p}{q} \right)} \frac{4}{\Gamma \left( \frac{2q+p}{q} \right)}.
\end{align*}
The exact solution can then be expressed as
\begin{align*}
x_{1, 2} (t) &= A \left( t-\frac{1}{3} \right)^{2 \nu} + B \left( t-\frac{1}{3} \right)^{2 \nu +1}, \\
x_{2, 2} (t) &= C \! + \! D \left( t -\frac{1}{3} \right)^{\nu} \! + E \left( t-\frac{1}{3} \right)^{\nu+1} \! + F \left( t-\frac{1}{3} \right)^{2 \nu} \\
& \ + G \left( t-\frac{1}{3} \right)^{2 \nu +1} + H \left( t-\frac{1}{3} \right)^{2 \nu +2}.
\end{align*}
Following the algorithm, exact solution of the problem can be found on further intervals.

For $\nu =1$, the solution on the interval $\displaystyle{\left[ 0, \frac{1}{3} \right]}$ is
\begin{align*}
x_{1, 1} (t) &= 0, \\
x_{2, 1} (t) &= t + t^2,
\end{align*}
while on the interval $\displaystyle{\left[ \frac{1}{3}, \frac{2}{3} \right]}$ we have
\begin{align*}
x_{1, 2} (t) &= \frac{1}{2} \left( t-\frac{1}{3} \right)^2 + \frac{1}{3} \left( t-\frac{1}{3} \right)^3, \\
x_{2, 2} (t) &= \frac{4}{9} \! + \! \frac{5}{3} \left( t -\frac{1}{3} \right) \! + \frac{4}{3} \left( t-\frac{1}{3} \right)^2  + \frac{8}{9} \left( t-\frac{1}{3} \right)^3 \\
& \ + \frac{1}{2} \left( t-\frac{1}{3} \right)^4.
\end{align*}
Expanding $x_{1, 2}$ and $x_{2, 2}$ we get
\begin{align*}
x_{1, 2} (t) &= \frac{7}{162} - \frac{2}{9} t+ \frac{1}{6} t^2 + \frac{1}{3} t^3, \\
x_{2, 2} (t) &= \frac{5}{486} + t + \frac{7}{9} t^2 + \frac{2}{9} t^3 + \frac{1}{2} t^4.
\end{align*}
We can see that we obtained the same exact solution as the solution presented in the paper by Rahimkhani et al. \cite{rahimkhani}.

\section{Conclusion}\label{conclusion}
Numerical solution of fractional control problems with constant delays in the state was discussed in the paper. An algorithm for fractional systems of single rational order with commensurate delays was introduced. The fractional system is transformed into a system of recurrence relations which can be easily solved by computer. The approximate solution is a truncated power series in each interval of the method of steps.

An example was provided to demonstrate that the results of the paper are convenient to solve concrete problems with given parameters. Numerical solution of a two-dimensional single-order system with two commensurate delays in the state variable and arbitrary order $\nu$ between $0$ and $1$ was found on the interval $[0, \frac{2}{3}]$. Since both the initial complete state and the control are polynomials, our approximate solution coincides with the exact solution of the problem. When $\nu$ approaches $1$, the numerical solution of the fractional-order system with Caputo derivative converges to the solution of integer-order system of order $1$ as expected. It means that the presented method produces reliable results which are in a good concordance with results known for integer-order systems.

The algorithm can be further generalized to provide a numerical scheme for solving fractional control problems with delays in control or both state and control, as well as for multi-order systems or equivalent multi-term equations. An open question is to find a generalization of the presented method for fractional systems with non-constant delays and distributed delays.


\section*{Acknowledgment}

The research was supported by the Czech Science Foundation under the project 16-08549S. This support is gratefully acknowledged.



%\bibliographystyle{IEEEtran}

\bibliographystyle{IEEEbib/IEEEtran}
\bibliography{EECS-CST2017_Rebenda_Smarda}

\end{document}